\documentclass[10pt]{amsart}
\usepackage{amsmath}
\usepackage{amsxtra}
\usepackage{amscd}
\usepackage{amsthm}
\usepackage{amsfonts}
\usepackage{amssymb}
\usepackage{eucal}

\newtheorem{conj}[subsection]{Conjecture}

\newtheorem{thm}[subsection]{Theorem}

\newcommand{\conjref}[1]{Conjecture~\ref{#1}}

\newcommand{\thmref}[1]{Theorem~\ref{#1}}
\newcommand{\secref}[1]{Sect~\ref{#1}}

\newcommand{\nc}{\newcommand}

\nc{\ssec}{\subsection}

\nc{\renc}{\renewcommand}
\nc{\on}{\operatorname}
\nc\ol{\overline}
\nc\wt{\widetilde}
\nc{\Loc}{\on{Loc}}
\nc{\Bun}{\on{Bun}}
\nc{\BQ}{{\mathbb{Q}}}
\nc{\BC}{{\mathbb{C}}}
\nc{\BG}{{\mathbb{G}}}
\nc{\CA}{{\mathcal{A}}}
\nc{\CC}{{\mathcal{C}}}
\nc{\CO}{{\mathcal{O}}}
\nc{\CR}{{\mathcal{R}}}
\nc{\CW}{{\mathcal{W}}}
\nc{\D}{{\mathcal{D}}}
\nc{\fg}{{\mathfrak{g}}}
\nc{\fD}{{\mathfrak{D}}}
\nc{\fh}{{\mathfrak{h}}}
\nc{\fn}{{\mathfrak{n}}}
\nc{\sM}{{\mathsf M}}
\nc{\ppart}{(\!(t)\!)}
\nc{\hg}{{\widehat\fg}}

\nc{\fW}{{\mathfrak{W}}}
\nc{\reg}{{\text{\rm reg}}}
\nc{\nilp}{{\text{\rm nilp}}}
\nc{\cG}{{\check{G}}}
\nc{\cB}{{\check{B}}}
\nc{\cg}{{\check{\fg}}}
\nc{\cb}{{\check{\fb}}}
\nc{\cn}{{\check{\fn}}}
\nc{\mer}{{\on{mer}}}
\nc{\Const}{\mathsf{Const}}
\nc{\Whit}{\on{Whit}}
\nc{\KL}{\on{KL}}
\nc{\FS}{\on{FS}}
\nc{\LocSys}{\on{LocSys}}
\nc{\QCoh}{\on{QCoh}}
\nc{\Cat}{\on{Cat}}
\nc{\Op}{\on{Op}}
\nc{\Gr}{\on{Gr}}
\nc{\Fl}{\on{Fl}}
\nc{\Rep}{\on{Rep}}
\renc{\mod}{{\on{-mod}}}
\nc{\Conn}{\on{Conn}}

\begin{document}

\title{Quantum Langlands Correspondence}

\author{Dennis Gaitsgory}

\date{Nov 13, 2007}

\maketitle

\bigskip

This note summarizes some conjectures on the theme of quantum 
geometric Langlands correspondence, which arose in the course
of discussions around October-November 2006 between A.~Beilinson, 
R.~Bezrukavnikov, A.~Braverman, M.~Finkelberg, D.~Gaitsgory and J.~Lurie.
We also thank E.~Frenkel and E.~Witten for stimulating
conversations.
\footnote{Tho sole bearer of responsibility for the speculations
that follow is D.G., who took it upon himself to write them down.
Other people mentioned are responsible only insofar as they
personally choose to do so. Accordingly, we shall not specify
individual credits, except for Jacob's initial guess of how to get the 
quantum group from twisted Whittaker sheaves.}

\medskip

Let $G$ be a reductive group over a field of characteristic $0$, and let
$\cG$ be its Langlands dual. By a level $c$ we will mean a choice of
a symmetric invariant form on the Lie algebra $\fg$. We will absorb
the critical shift into the notation, i.e., $c=0$ means the critical level.
Given $c$, we will denote by $\frac{1}{c}$ the dual level for $\cG$,
obtained by identifying the corresponds Cartan subalgebras as
duals one one another.

\section{Categories acted on by the loop group}

We will assume having the following notions at our
disposal:

\ssec{} The notion of category $\CC$ acted on by the loop group $G\ppart$
at level $c$ (for an abelian category this notion is developed, e.g., in \cite{FG2}).

\ssec{Example} The category $\fD^c(G\ppart)\mod$ of ($c$-twisted) D-modules
on the loop group itself, carries an action of $G\ppart$ at level $c$ and
a commuting action at level $-c$. From now on we will denote it
by $\fD^{c,-c}(G\ppart)\mod$. We have a canonical equivalence
$$\fD^{c,-c}(G\ppart)\mod\simeq \fD^{-c,c}(G\ppart)\mod,$$
that interchanges the two actions.

\ssec{} To $\CC_1$ acted on by $G\ppart$ at level $c$, and $\CC_2$
acted on by $G\ppart$ at level $-c$, we should be able to associate their
tensor product over $G\ppart$, denoted $\CC_1\underset{G\ppart}\otimes \CC_2$.

\ssec{Example} let $K$ be a subgroup of $G[[t]]$,
and let $\CC_1=\fD^{c}(G\ppart/K)\mod$. Then for $\CC_2:=\CC$ as above,
$\fD^{c}(G\ppart/K)\mod \underset{G\ppart}\otimes \CC$ should be equivalent to the
category $\CC^K$ of $K$-equivariant objects in $\CC$. In particular, we have:
\begin{equation} \label{D unit}
\fD^{c,-c}(G\ppart)\mod \underset{G\ppart}\otimes \CC\simeq \CC,
\end{equation}
as categories acted on by $G\ppart$ at level $-c$.

\ssec{Example} \label{weak eq}
Let $\CC_1=\hg^c\mod$--the category of smooth Kac-Moody
modules at level $c$. For $\CC_2:=\CC$ consider the category 
$\hg^c\mod \underset{G\ppart}\otimes \CC$.

\medskip

If instead of $G\ppart$ we had a 
group-scheme $H$, the corresponding category $\fh\mod\underset{H}\otimes \CC$
would identify with the category $\CC^{H,w}$ of weakly $H$-equivariant objects in $\CC$.
In this case, the tensor category $\Rep(H)$ would act on $\CC^{H,w}$. 
In the sequel, we will see what replaces this structure when instead of
$H$ we have the loop group $G\ppart$.

\ssec{} To $\CC$ acted on by $G\ppart$ at level $c$, we should be able to assign 
the category $\Whit(\CC)$ that corresponds to $(N((t)),\chi)$-equivariant objects, where
$\chi:N((t))\to \BG_a$ is a non-degenerate character.

\ssec{Example} One of the principal players for us will be the category
\begin{equation} \label{Whit on Gr}
\Whit^c(\Gr_G):=\Whit(\fD^{c}(\Gr_G)\mod).
\end{equation}
Another important example is
$\Whit^c(G\ppart):=\Whit(\fD^{c,-c}(G\ppart)\mod)$. By transport of structure,
the latter carries an action of $G\ppart$ at level $-c$. We have:
\begin{equation} \label{Whit as tensor}
\Whit(\fD^{c,-c}(G\ppart)\mod)\underset{G\ppart}\otimes \CC\simeq \Whit(\CC).
\end{equation}

\ssec{Example} Consider $\CC=\hg^c\mod$. We are supposed to have
\begin{equation} \label{W}
\Whit(\hg^c\mod)\simeq \CW^c_\fg\mod,
\end{equation}
where $\CW^c_\fg$ is the W-algebra corresponding to $\fg$ at level $c$.

\section{Chiral categories}

\ssec{} Let $X$ be an algebraic curve. Another notion that we assume having at 
our disposal is that of chiral category over $X$. 

\medskip

The data of a chiral category $\CA$ assigns to each integer $n$ an $\CO$-module of categories
$\CA^n$ over $X^n$, equipped with factorization isomorphisms, which we will spell out for $n=2$:

\medskip

\noindent The restriction $\CA^2|_{X\times X-\Delta(X)}$ should be identified
with $\CA^1\boxtimes \CA^1|_{X\times X-\Delta(X)}$, and the restriction $\CA^2|_{\Delta(X)}$
should be identified with $\CA^1$. 

\medskip

In addition, $\CA^n$ is supposed to have a unit object, analogously to the case of
chiral algebras. The latter endows the sheaf of categories $\CA^n$ with a connection along $X$. 

\ssec{} We will usually think of a chiral category as an $\CO$-module of categories $\CA^1$
over $X$ itself, endowed with an extra structure. When a chiral category is obtained by a 
universal procedure (i.e., is a vertex operator category), we will think of it as a plain
category $\CA$ equal to the fiber of $\CA^1$ at a point $x\in X$, 
endowed with an extra structure.

\medskip

The notion of chiral category should be regarded as a D-module version of the notion of 
$E_2$-category. 

\ssec{Example} A sheaf of symmetric monoidal categories over $X$, endowed with a
connection along $X$ gives rise to a chiral category. 

\ssec{} When working over $\BC$, and we choose a coordinate on our curve,
there is a transcendental procedure that assigns to a ribbon category a chiral 
category. (We are not going to use it.) This procedure is fully faithful: a chiral
category comes in this was from a ribbon category when a certain representability
condition is satisfied. This is how a monoidal structure arises in \cite{KL}.

\ssec{Example} \label{chiral alg}
Let $A$ be a chiral algebra. Then the category $A\mod$
of chiral $A$-modules is naturally a chiral category. 

\ssec{Example} A construction of Beilinson and Drinfeld endows the 
category $\fD^c(\Gr_{G}\mod)$ with a structure of chiral category.

\ssec{Example} The category $\Whit^c(\Gr_G)$ is a chiral category.

\ssec{} By analogy with the theory of chiral algebras, given a chiral category, it makes
sense to consider module categories over it. We will consider module categories
supported at a fixed point $x$ of the curve, with $t$ being a local coordinate.

\ssec{} Generalizing the construction of Beilinson and Drinfeld, we obtain that
the category $\fD^c(G\ppart\mod)$ is naturally a module category with respect to
$\fD^{c}(\Gr_{G}\mod)$. In addition, it carries a commuting action of $G\ppart$ at level $-c$.
Hence, for any category $\CC$ acted on by $G\ppart$ at level $c$, the category
$\fD^{c,-c}(G\ppart)\mod \underset{G\ppart}\otimes \CC$ is a module category
with respect to $\fD^c(\Gr_{G}\mod)$. Taking into account \eqref{D unit}, we obtain
that on a given category $\CC$ an action of $G\ppart$ at level $c$ gives rise to
a structure of module category with respect to $\fD^c(\Gr_{G}\mod)$.

\medskip

\begin{conj} For a category $\CC$, the data of an action of $G\ppart$ at level $c$ is
equivalent to a structure of module category with respect to $\fD^c(\Gr_{G}\mod)$.
\end{conj}

\ssec{Example} Independent of the above conjecture, the category $\Whit^c(G\ppart)$
has a structure of module category with respect to $\Whit^c(\Gr_G)$, and carries
a commuting action of $G\ppart$ at level $-c$. 

\medskip

Hence, by \eqref{Whit as tensor}, for
any category $\CC$ acted on by $G\ppart$ at level $c$, the category $\Whit(\CC)$ is a module
category with respect to $\Whit^c(\Gr_G)$.

\ssec{} By \secref{chiral alg}, the category $\hg^c\mod$ is a chiral category. 
Let $\KL_G^c\subset \hg^c\mod$ be the subcategory
consisting of $G[[t]]$-integrable representations. (The symbol $\KL$ stands
for Kazhdan-Lusztig, who studied this category in \cite{KL}.) I.e.,
$$\KL_G^c:=(\hg^c\mod)^{G[[t]]}.$$
This is also a chiral category.
We will regard the fiber of $\hg^c\mod$ at $x\in X$ as a module category
with respect to $\KL^c_G$.

\medskip

The following is established in \cite{FG1}:

\begin{thm}  \label{FG3}
The following two pieces of structure defined on $\hg^c\mod$ commute:
the action of $G\ppart$ at level $c$, and the structure of module category
with respect to $\KL_G^c$.
\end{thm}

\ssec{} Let $\CC$ be again a category acted on by $G\ppart$ at level $-c$, and
consider the category $\KL(\CC):=\hg^c\mod \underset{G\ppart}\otimes \CC$ of
\secref{weak eq}. From \thmref{FG3} we obtain that this category is
naturally a module category with respect to $\KL_G^c$.

\medskip

It is this structure that we regard as a substitute for the action of $\Rep(H)$,
alluded to in \secref{weak eq}.

\section{Local Quantum Langlands}

\ssec{}

The following is a version of a conjecture proposed by J.~Lurie:

\begin{conj}  \label{Lurie}
For $c$ not rational negative, the chiral categories $\Whit^c(\Gr_G)$ and
$\KL^{\frac{1}{c}}_\cG$ are equivalent.
\end{conj}

This conjecture has been essentially proven in \cite{Ga} for $c$ irrational,
by identifying both sides with a third chiral category, namely, that of factorizable
sheaves of \cite{BFS}. 

\ssec{}

We can now formulate the Local Quantum Geometric Langlands conjecture,
which we literally believe to hold only for irrational values of $c$:

\begin{conj}  \label{local Langlands}
There exists a 2-equivalence $\Psi^{c,-\frac{1}{c}}_{G\to \cG}$:
$$\{\text{ Categories acted on by $G\ppart$ at level $c$ }\} \to 
\{\text{ Categories acted on by $\cG\ppart$ at level $-\frac{1}{c}$ }\},$$
characterized by {\bf either} of the following two properties: for $\CC$ acted on by $G\ppart$
and $\check\CC:=\Psi^{G\to \cG}(\CC)$, we need that:

\begin{itemize}

\item The category $\Whit(\CC)$, regarded as
a module category with respect to $\Whit^c(\Gr_G)$, is equivalent to
$\KL(\check\CC)$, regarded as a module category
with respect to $\KL_\cG^{\frac{1}{c}}$, when we identify 
$\Whit^c(\Gr_G)\simeq \KL_\cG^{\frac{1}{c}}$ via \conjref{Lurie}.

\item The category $\KL(\CC)$, regarded as
a module category with respect to $\KL_G^{-c}$, is
equivalent to $\Whit(\check\CC)$, regarded as a module category
with respect to $\Whit^{-\frac{1}{c}}(\Gr_\cG)$, when we identify 
$\KL_G^{-c}\simeq \Whit^{-\frac{1}{c}}(\Gr_\cG)$ via \conjref{Lurie}.

\end{itemize}

\end{conj}

\ssec{}

We can divide \conjref{local Langlands} into three steps:

\begin{conj} \label{Whit conj}
For $c$ irrational, the assignment $\CC\mapsto \Whit(\CC)$ establishes a 2-equivalence
$$\{\text{Categories acted on by $G\ppart$ at level $c$}\} \to 
\{\text{Module categories with respect to $\Whit^c(\Gr_G)$}\}.$$
\end{conj}

\begin{conj} \label{KL conj}
The assignment $\CC\mapsto \KL(\CC)$ establishes a 2-equivalence
$$\{\text{Categories acted on by $G\ppart$ at level $c$}\} \to 
\{\text{Module categories with respect to $\KL^{-c}_G$}\}.$$
\end{conj}

\begin{conj} \label{symmetry}
Assuming the above two conjectures, the two composed 2-functors
$$\left(\KL \text{ for $\cG$ }\right)^{-1} \circ \left(\Whit\text{ for $G$ }\right)$$
and
$$\left(\Whit \text{ for $\cG$ }\right)^{-1} \circ \left(\KL\text{ for $G$ }\right)$$
$$\{\text{ Categories acted on by $G\ppart$ at level $c$ }\} \to 
\{\text{ Categories acted on by $\cG\ppart$ at level $-\frac{1}{c}$ }\},$$
are isomorphic.
\end{conj}

\ssec{}

Let us consider some examples of how the 2-functor $\Psi^{c,-\frac{1}{c}}_{G\to \cG}$ is
supposed to act. We claim that we have:
$$\Psi^{c,-\frac{1}{c}}_{G\to \cG}(\fD^c(\Gr_G)\mod)\simeq \fD^{-\frac{1}{c}}(\Gr_\cG)\mod.$$
This follows from either of the characterizing properties of $\Psi$, since
$$\Whit(\fD^c(\Gr_G)\mod):=\Whit^c(\Gr_G)$$ and 
$$\KL(\fD^c(\Gr_G)\mod):=\hg^{-c}\mod \underset{G\ppart}\otimes \fD^c(\Gr_G)\mod
\simeq \hg^{-c}\mod^{G[[t]]}=:\KL_G^{-c}.$$

\begin{conj}  \label{Iwahori conj}
$\Psi^{c,-\frac{1}{c}}_{G\to \cG}(\fD^c(\Fl_G)\mod)\simeq 
\fD^{-\frac{1}{c}}(\Fl_\cG)\mod$.
\end{conj}

This is an interesting conjecture in its own right, as it translates as:
\begin{conj} \label{Iwahori conj reform}
$$\Whit^{c}(\Fl_G)\simeq 
\widehat\cg^{\frac{1}{c}}\mod^{\check{I}}.$$
\end{conj}

\ssec{Duality of $\CW$-algebras}

Let us recall the assertion of \cite{FF} that there exists an isomorphism
\begin{equation} \label{W isom}
\CW_\fg^c\simeq \CW_\cg^{\frac{1}{c}}.
\end{equation}
In particular, the corresponding categories of modules are equivalent.

\medskip

Note, however, that the identification \eqref{W}
defines on $\CW_\fg^c\mod$ the commuting structures of module over the
chiral categories $\Whit^c(\Gr_G)$ and $\KL^c_G$.

\medskip

We propose the following strengthening of \eqref{W isom}:

\begin{conj}  \label{strong W}
The equivalence of categories
$$\CW_\fg^c\mod\simeq \CW_\cg^{\frac{1}{c}}\mod,$$
induced by \eqref{W isom}, respects the module structures with respect to
the chiral categories
$$\Whit^c(\Gr_G)\simeq \KL_\cG^{\frac{1}{c}} \text{ and }
\KL^c_G \simeq \Whit^{\frac{1}{c}}(\Gr_\cG).$$
\end{conj}

This conjecture formally implies that
\begin{equation} \label{Psi of Whit}
\Psi^{c,-\frac{1}{c}}_{G\to \cG}(\hg^{c}\mod)\simeq \Whit^{\frac{1}{c}}(\cG\ppart) 
\text{ and }
\Psi^{c,-\frac{1}{c}}_{G\to \cG}(\Whit^{c}(G\ppart))\simeq \widehat\cg\mod^{\frac{1}{c}}.
\end{equation}

\ssec{}

Let us denote by $\sM^{c,-\frac{1}{c}}_{G\to \cG}$ the category
$\Psi^{c,-\frac{1}{c}}_{G\to \cG}(\fD^{c,-c}(G\ppart)\mod)$. It carries an action
of $\cG$ at level $-\frac{1}{c}$ and a commuting action of $G\ppart$ at
level $-c$, and the functor $\Psi^{c,-\frac{1}{c}}$ can be realized as
$$\CC\mapsto \sM^{c,-\frac{1}{c}}_{G\to \cG}\underset{G\ppart}\otimes \CC.$$

\medskip

Isomorphism \eqref{Psi of Whit} implies that
\begin{equation} \label{extra symmetry}
\sM^{c,-\frac{1}{c}}_{G\to \cG}\simeq \sM^{\frac{1}{c},-c}_{\cG\to G}.
\end{equation}

\medskip

The category $\sM^{c,-\frac{1}{c}}_{G\to \cG}$ has the following properties
\begin{equation} \label{Whit to g}
\left(\Whit\text{ for $G$ }\right)(\sM^{c,-\frac{1}{c}}_{G\to \cG})\simeq \widehat\cg\mod^{-\frac{1}{c}}
\text{ and } 
\left(\Whit\text{ for $\cG$ }\right)(\sM^{c,-\frac{1}{c}}_{G\to \cG})\simeq \hg\mod^{-c}.
\end{equation}

\medskip

\noindent{\it Remark.} In \cite{Sto} it is suggested that the should exist a chiral algebra
$M^{c,-\frac{1}{c}}_{G\to \cG}$, which receives maps with commuting images
from the Kac-Moody chiral algebras $A_{\fg,-c}$ and $A_{\cg,-\frac{1}{c}}$,
corresponding to $\fg$ and $\cg$ at levels $-c$ and $-\frac{1}{c}$, respectively, 
such that the Drinfeld-Sokolov reduction of $M^{c,-\frac{1}{c}}_{G\to \cG}$ with respect
to $\fg$ is isomorphic to $A_{\cg,-\frac{1}{c}}$, and with respect
to $\cg$ is isomorphic to $A_{\fg,-c}$. The above discussion does not produce such
a chiral algebra, but rather a chiral category with the corresponding properties. 

\ssec{}

Equation \eqref{extra symmetry} implies that for $\CC_1$, acted on by $G\ppart$ at level
$c$, $\CC_2$, acted on by $G\ppart$ at level $-c$, and 
$$\check\CC_1:=\Psi^{c,-\frac{1}{c}}_{G\to \cG}(\CC_1),\,\,\,
\check\CC_2:=\Psi^{-c,\frac{1}{c}}_{G\to \cG}(\CC_2),$$
we have:
\begin{equation} \label{compat with ten prod}
\CC_1\underset{G\ppart}\otimes \CC_2\simeq 
\check\CC_1\underset{\cG\ppart}\otimes \check\CC_2.
\end{equation}

In particular, for $\CC$ acted on by $G\ppart$ at level $c$ and 
$\check\CC:=\Psi^{c,-\frac{1}{c}}_{G\to \cG}(\CC)$ we obtain:
\begin{equation} \label{sph preserved}
\CC^{G[[t]]}\simeq \CC\underset{G\ppart}\otimes \fD^{-c}(\Gr_G)\mod\simeq
\check\CC\underset{G\ppart}\otimes \fD^{\frac{1}{c}}(\Gr_\cG)\mod\simeq \check\CC^{\cG[[t]]}.
\end{equation}

Assuming \eqref{Iwahori conj}, we also obtain
\begin{equation} \label{Iw preserved}
\CC^I\simeq \check\CC^{\check{I}}.
\end{equation}

\ssec{Harish-Chandra bimodules}

Let us denote by $\on{HCh}_G^{c,-c}$ the category
$$\hg^{c}\mod\underset{G\ppart}\otimes \hg^{-c}\mod.$$

\medskip

We remark that for a group scheme $H$, the corresponding category
$\fh\mod\underset{H}\otimes \fh\mod$ is indeed tautologically equivalent
to the category of Harish-Chandra modules for the pair $(\fh\oplus\fh,H)$.

\medskip

Equation \eqref{Whit to g} implies that we have the equivalence
\begin{conj} \label{HCh}
$$\on{HCh}_G^{c,-c}\simeq (\Whit\times \Whit)
(\fD^{\frac{1}{c},-\frac{1}{c}}(\cG\ppart)\mod).$$
\end{conj}


\section{Global quantum Langlands}

\ssec{}

Assume now that $X$ is a complete curve. Let $\Bun_G$ denote the moduli
stack of $G$-bundles on $X$.

The following conjecture was
proposed in \cite{Sto}:

\begin{conj} \label{Sto}
There exists an equivalence of categories
$$\fD^c(\Bun_G)\mod \simeq \fD^{-\frac{1}{c}}(\Bun_\cG)\mod.$$
\end{conj}

We will now couple this with \conjref{Lurie}, which would, conjecturally,
fix the equivalence of \conjref{Sto} uniquely.

\ssec{}

Let $x_1,...,x_n$ be a finite collection of points on $X$. On the one hand, 
we have a localization functor
$$\on{Loc}:\KL^c_{G,x_1}\times...\times \KL^c_{G,x_n}\to \fD^c(\Bun_G),$$
obtained by considering conformal blocks of $\hg^c$-modules.

\medskip

On the other hand, we the Poincare series functor
$$\on{Poinc}:\Whit^{-c}(\Gr_{G,x_1})\times...\times \Whit^{-c}(\Gr_{G,x_n})\to \fD^c(\Bun_G),$$
corresponding to the diagram
$$\Gr_{G,x_1}/N\ppart \times...\times \Gr_{G,x_n}/N\ppart \leftarrow
\Gr_{G,x_1}\times...\times \Gr_{G,x_n}/N_{out}\to 
\Gr_{G,x_1}\times...\times \Gr_{G,x_n}/G_{out},$$
where
$$(\Gr_{G,x_1}\times...\times (\Gr_{G,x_n})/G_{out}
\simeq \Bun_G.$$

\ssec{Global unramified quantum Langlands}

We propose:

\begin{conj}  \label{global unramified quantum Langlands}
There exists an equivalence as in \conjref{Sto}, which for every
collection of points $x_1,...,x_n$ makes the diagram
$$
\CD
\KL^c_{G,x_1}\times...\times \KL^c_{G,x_n}  @>{\text{\conjref{Lurie}}}>> 
\Whit^{\frac{1}{c}}(\Gr_{\cG,x_1})\times...\times \Whit^{\frac{1}{c}}(\Gr_{\cG,x_n}) \\
@V{\on{Loc}}VV    @V{\on{Poinc}}VV  \\
\fD^c(\Bun_G)\mod @>{\text{\conjref{Sto}}}>>   \fD^{-\frac{1}{c}}(\Bun_\cG)\mod
\endCD
$$
commute.
\end{conj}

\ssec{The ramified case}

For a point $x\in X$ let $\Bun_{G,x}$ be the moduli space of $G$-bundles with
a full level structure at $x$. The category $\fD^c(\Bun_{G,x})\mod$ carries a natural
action of $G\ppart$ at level $c$.  

\medskip

We propose the following:

\begin{conj} \label{global ramified quantum Langlands}
The 2-functor $\Psi^{c,-\frac{1}{c}}_{G\to \cG}$ sends the
category $\fD^c(\Bun_{G,x})\mod$ to $\fD^{-\frac{1}{c}}(\Bun_{\cG,x})$.
\end{conj}

We emphasize that this conjecture {\it does not} say that the categories
$\fD^c(\Bun_{G,x})\mod$ and $\fD^{-\frac{1}{c}}(\Bun_{\cG,x})$ are equivalent.
Rather, it says that they correspond to each other under $\Psi^{c,-\frac{1}{c}}_{G\to \cG}$.

\ssec{Compatibility with the unramified picture}

Let us couple \conjref{global ramified quantum Langlands}
with \eqref{sph preserved} and \eqref{Iw preserved}.
We obtain
\begin{equation} \label{global sph}
\left(\fD^c(\Bun_{G,x})\mod\right)^{G[[t]]}\simeq 
\left(\fD^{-\frac{1}{c}}(\Bun_{\cG,x})\right)^{\cG[[t]]}
\end{equation}
and
\begin{equation} \label{global Iw}
\left(\fD^c(\Bun_{G,x})\mod\right)^{I}\simeq 
\left(\fD^{-\frac{1}{c}}(\Bun_{\cG,x})\right)^{\check{I}},
\end{equation}
respectively.

\medskip

However, $\left(\fD^c(\Bun_{G,x})\mod\right)^{G[[t]]}\simeq \fD^c(\Bun_{G})\mod$,
so \eqref{global sph} recovers \conjref{Sto}.

\medskip

Note that $\left(\fD^c(\Bun_{G,x})\mod\right)^{I}\simeq \fD^c(\Bun'_{G})\mod,$
where $\Bun'_{G}$ denotes the moduli space of $G$-bundles with a parabolic
structure at $x$. So, \eqref{global Iw} leads to an equivalence
$$\fD^c(\Bun'_{G})\mod\simeq \fD^{-\frac{1}{c}}(\Bun'_{\cG})\mod,$$
generalizing \eqref{Sto}.

\ssec{}

One can give a version of \conjref{global ramified quantum Langlands} along
the lines of \conjref{global unramified quantum Langlands}. Let us instead
of one point $x$ have a collection $x_1,...,x_n$. We can consider the functors
$$\on{Loc}_{ram}:\hg_{x_1}^c\mod\times...\times \hg_{x_n}^c\mod\to \fD^c(\Bun_{G,x_1,...,x_n})\mod$$
and 
$$\on{Poinc}_{ram}:\Whit^{-c}(G\ppart_{x_1})\times...\times \Whit^{-c}(G\ppart_{x_n})\to 
\fD^c(\Bun_{G,x_1,...,x_n})\mod.$$

\begin{conj}
The 2-functor $\Psi^{c,-\frac{1}{c}}_{G\to \cG}$ takes the functor $\on{Loc}_{ram}$ for
$G$ at level $c$ to the functor $\on{Poinc}_{ram}$ for $\cG$ at level $-\frac{1}{c}$.
\end{conj}

\section{Local correspondence at $c=0/\infty$}

Let us now consider the limiting cases, that correspond to the
"classical", i.e., non-quantum, geometric Langlands correspondence.
We remind that $c=0$ means the critical level.

\ssec{}

First, some comments are due as to how the corresponding objects
look at $c=\infty$.

\medskip

Let $\Conn_G(\D^\times)$ denote the ind-scheme of $G$-connections over the formal
punctured disc $\D^\times$. We have an action of $G\ppart$ on $\Conn(\D^\times)$
by gauge transformations. 

\medskip

By definition, a category $\CC$ acted on by $G\ppart$ at level $\infty$ is a category
over $\Conn_G(\D^\times)$, equipped with a compatible {\it weak} action of $G\ppart$.

\medskip

The category $\hg^\infty\mod$ is by definition $\QCoh(\Conn_G(\D^\times))$. 

\medskip

The category
$\KL^\infty_G$ is the category of quasi-coherent sheaves on $\Conn_G(\D^\times)$ that
are supported on subscheme $\Conn^\reg_G(\D^\times)$ of regular connections
(=without poles) and that are equivariant with respect to $G[[t]]$. I.e., this is the
category $\QCoh(\Conn^\reg_G(\D^\times)/G[[t]])$. However, 
$\Conn^\reg_G(\D^\times)/G[[t]]\simeq \on{pt}/G$, so 
$$\KL^\infty_G\simeq \Rep(G).$$

Similarly, 
$$\hg^\infty\mod^I\simeq \QCoh(\fn/B).$$



\ssec{}

\conjref{KL conj} translates into the following:

\begin{conj} \label{KL conj infty}
The assignment $\CC\mapsto \CC^{G\ppart,w}$ establishes
a 2-equivalence between the 2-category of categories over $\Conn_G(\D^\times)$
equipped with a weak $G\ppart$ action, and the 2-category of modules
with respect to the chiral category $\Rep(G)$.
\end{conj}

Let $\LocSys_G(\D^\times)$ denote the quotient stack $\Conn_G(\D^\times)/G\ppart$.

\medskip

Note that the assignment $\CC\mapsto \CC^{G\ppart,w}$ can be
alternatively interpreted as a bijection between the 2-category of
categories over $\Conn_G(\D^\times)$ equipped with a weak $G\ppart$ 
action, and the 2-category of categories over the stack 
$\LocSys_G(\D^\times)$.

\medskip

Hence, \conjref{KL conj infty} implies:

\begin{conj}  \label{interp over LocSys}
For a category $\CC$ the following two pieces of structure are equivalent:
a structure of category over the stack $\LocSys_G(\D^\times)$, and
a structure of module category with respect to $\Rep(\cG)$.
\end{conj}

\ssec{}
The category $\CW^\infty_G\mod$ identifies with $\QCoh(\Op_G(\D^\times))$, where 
$\Op_G(\D^\times)$ is the ind-scheme of $G$-opers on $\D^\times$.

\medskip

The category $\Whit^\infty(\Gr_G)$ identifies with $\QCoh(\Op^{unr}_G(\D^\times))$,
where $\Op^{unr}_G(\D^\times)$ is the ind-scheme opers that are unramified
as local systems.

\medskip

Similarly, the category $\Whit^\infty(\Fl_G)$ identifies with $\QCoh(\Op^{nlp.ram.}_G(\D^\times))$,
where we denote by $\Op^{nilp.ram.}_G(\D^\times)$ is the ind-scheme of opers that
have a nilpotent ramification as local systems.

\medskip

\noindent{\it Remark.} The ind-schemes $\Op^{unr}_G(\D^\times)$ and 
$\Op^{nilp.ram.}_G(\D^\times)$ should not be confused with their sub-schemes
$\Op^\reg_G$ and $\Op^\nilp_G$ that correspond to opers with a regular
and nilpotent singularity, respectively.

\medskip

Finally, the category $(\Whit\times \Whit)(\fD^{\infty,-\infty}(G\ppart)\mod)$ identifies with
$\QCoh(\on{Isom}_G(\D^\times))$, where $\on{Isom}_G(\D^\times)$ is the isomonodromy groupoid over 
$\Op_G(\D^\times)$.

\ssec{}

Let us now specialize some of the conjectures mentioned above
to $c=0$ and $\infty$.

\medskip

\conjref{Lurie} for $c=0$ reads
$$\Whit^0(\Gr_G)\simeq \Rep(\cG),$$
as chiral categories. This is a valid assertion.

\medskip

\conjref{Lurie} for $c=\infty$ reads
$$\KL^0_G\simeq \QCoh(\Op^{unr}_\cG(\D^\times)).$$
This is also a theorem, established in \cite{FG3}.

\medskip

\conjref{Iwahori conj reform} for $c=0$ reads as
$$\Whit^0(\Fl_G)\simeq \QCoh(\cn/\cB),$$
which is the theorem of \cite{AB}.

\medskip

\conjref{Iwahori conj reform} for $c=\infty$ reads as
\begin{conj}
$$\hg^0\mod^I\simeq \QCoh(\Op^{nlp.ram.}_\cG(\D^\times)).$$
\end{conj}
The latter can be viewed as a generalization of the
main conjecture from \cite{FG2}.

\medskip

\conjref{HCh} for $c=0$ reads
\begin{conj}
$$\hg^0\mod\underset{G\ppart}\otimes \hg^0\mod\simeq
\QCoh(\on{Isom}_G(\D^\times)).$$
\end{conj}

\ssec{}

Consider now the 2-functor $\Psi^{0,\infty}_{G\to \cG}$. We will compose
it with the 2-equivalence of the (plausible) \conjref{KL conj infty}
and \conjref{interp over LocSys},
and thus consider the 2-functor, denoted $\Phi_{G\to \cG}$ from
the 2-category of categories acted on by $G\ppart$ at level $0$
to that of categories over the stack $\LocSys_G(\D^\times)$.

\medskip

By construction, the 2-functor in question is
$$\CC\mapsto \Whit(\CC),$$
when the latter is regarded as a module category with respect to
the chiral category $\Rep(\cG)$.

\medskip

\noindent{\it Remark.}
At level $0$, the 2-functor $\Whit$ is not a 2-equivalence for obvious
reasons: it kills the category $\CC=\on{Vect}$, equipped with the trivial
action of $G\ppart$. Thus, in order to have a Langlands-type equivalence
in this case, one has to enhance the RHS, presumably by adding
an Arthur $SL_2$.

\section{Global correspondence at $c=0/\infty$}

\ssec{The unramified case}

First, we note that the equivalence \eqref{Sto} specializes at $c=0$
to the usual geometric Langlands:

\begin{conj}   \label{usual Langlands}
There exists an equivalence
$$\fD^0(\Bun_G)\mod \simeq \on{QCoh}(\LocSys_\cG(X)).$$
\end{conj}

\medskip

The commutativity of the diagram in 
\conjref{global unramified quantum Langlands} for $c=0$ amounts to the
Beilinson-Drinfeld construction of Hecke eigensheaves via opers.

\medskip

When we exchange the roles of $G$ and $\cG$ and replace $0$ by
$\infty$, the commutative diagram of 
\conjref{global unramified quantum Langlands} amounts to the
expectation that the equivalence of \eqref{usual Langlands} takes
the "Whittaker coefficient" D-modules on the LHS to the tautological
coherent sheaves associated to points of $X$ and representations
of $\cG$ on the RHS.

\ssec{The ramified case}

Consider the 2-functor $\Phi_{G\to \cG}$ applied to the category
$\fD^0(\Bun_{G,x})$. This is the category $\Whit(\fD^0(\Bun_{G,x}))$
over $\LocSys_G(\D^\times)$. By \conjref{global ramified quantum Langlands}
this category is equivalent to
$$\on{QCoh}(\LocSys_{\cG,x}(X))^{\cG\ppart,w}.$$

\medskip

Note that the stack $\LocSys_{\cG,x}(X)$ classifies $\cG$-local
systems on $X$ with an arbitrary ramification at $x$, and a 
full level structure at $x$ on the underlying $G$-bundle. Hence,
the stack $\LocSys_{\cG,x}(X)/\cG\ppart$ identifies with the
stack $\LocSys_\cG(X-x)$ of $\cG$-local systems defined
on the punctured curve. The category 
$$\on{QCoh}(\LocSys_{\cG,x}(X))^{\cG\ppart,w}\simeq
\on{QCoh}(\LocSys_{\cG}(X-x))$$
is naturally a category over $\LocSys_G(\D^\times)$ via
the map of stacks $$\LocSys_{\cG}(X-x)\to \LocSys_G(\D^\times).$$

\medskip

Summarizing, we obtain:

\begin{conj}
The category $\Whit(\fD^0(\Bun_{G,x}))$ is equivalent to
$\on{QCoh}(\LocSys_{\cG}(X-x))$, as categories over 
$\LocSys_G(\D^\times)$.
\end{conj}

\end{document}